\documentclass[a4paper,10pt]{amsart}

\usepackage{amssymb, amsthm, amsmath}
\usepackage{algorithmic}

\newtheorem{theorem}{Theorem}[section]
\newtheorem{lemma}[theorem]{Lemma}

\theoremstyle{definition}
\newtheorem{definition}[theorem]{Definition}
\newtheorem{example}[theorem]{Example}
\newtheorem{corollary}[theorem]{Corollary}

\theoremstyle{remark}
\newtheorem{remark}[theorem]{Remark}

\numberwithin{equation}{section}

\newcommand{\Z}{\mathbb{Z}}
\newcommand{\mL}{\mathcal{L}}

\DeclareMathOperator{\coeff}{coeff}
\DeclareMathOperator{\cond}{cond}

\begin{document}

\title{Tschirnhaus-Weierstrass curves}


\author{Josef Schicho}
\address{Johann Radon Institute for Computational and Applied Mathematics (RICAM)\\Austrian Academy of Sciences\\Altenbergerstrasse 69\\A-4040 Linz, Austria}
\email{josef.schicho@oeaw.ac.at}
\thanks{The first author is partially supported by the FWF (Austrian Science Fund) in the frame of project 18992.}

\author{David Sevilla}
\address{Johann Radon Institute for Computational and Applied Mathematics (RICAM)\\Austrian Academy of Sciences\\Altenbergerstrasse 69\\A-4040 Linz, Austria}
\email{david.sevilla@oeaw.ac.at}
\urladdr{http://www.davidsevilla.com}
\thanks{The second author is partially supported by the Spanish MEC project MTM2007-67088 and the FWF project P22766-N18.}


\date{}

\dedicatory{}

\begin{abstract}
We define the concept of Tschirnhaus-Weierstrass curve, named after the Weierstrass form of an elliptic curve and Tschirnhaus transformations. Every pointed curve has a Tschirnhaus-Weierstrass form, and this representation is unique up to a scaling of variables. This is useful for computing isomorphisms between curves.
\end{abstract}

\maketitle

\section{Introduction}

A elliptic curve over a characteristic zero field is given by a polynomial equation of degree 3 in two variables
\[\sum_{i+j\leq 3} a_{i,j} x^iy^j\]
that can be reduced, by invertible transformations, to the classical \emph{Weierstrass normal form of an elliptic curve}:
\[y^2=x^3+ax+b\]
which characterizes it.

Another definition of an elliptic curve is that of a genus 1 curve. Then, the Riemann-Roch theorem \cite{hartshorne77} determines the dimensions of the Riemann-Roch spaces, from which a certain equation of degree 3 is obtained. This equation can be simplified to obtain, again, the form above (we elaborate on this in Section \ref{SECTION: motivation and example}). It is of interest that imposing the previous form (that is, all the other coefficients of degree up to 3 are zero), the generators of the Riemann-Roch spaces are constrained in such a way that they are uniquely determined up to products by constants.

Our goal in this paper is to generalize the concept of Weierstrass nomal form for curves of higher genus, in such as a way as to also fix the generators of the Riemann-Roch spaces up to product by constants. The interest of such property lies in the fact that it allows us to improve an algorithm by F. Hess \cite{hess04} for computing the set of birational isomorphisms between two given algebraic curves. (If the two curves are equal, then the result is the automorphism group; if not, then the algorithm decides whether the two curves are birationally equivalent.)

Here is a summary of the main idea of Hess's algorithm. First, pick a Weierstrass point $P$ on the first curve (see Definition \ref{DEF: W place}). An isomorphism has to take $P$ to a Weierstrass point on the second curve with the same gap sequence. There are only finitely many such points, and we compute for each of them the set of isomorphisms mapping $P$ to it. These isomorphisms induce isomorphisms on corresponding Riemann-Roch spaces, and this fact allows one to compute the isomorphisms efficiently. It is surprising that the algorithm does not use the canonical embedding; the reason is that the embedding based on Weierstrass points leads to better complexity.

Hess's algorithm has been implemented in Magma \cite{magma}, and the algorithm is quite useful for practical computations when the coefficient field is finite. Our improvement applies in the characteristic zero case. We show that the birational isomorphisms of two Tschirnhaus-Weierstrass curves (see Section \ref{SECTION: TW curves}) mapping the unique point at infinity to the unique point at infinity can always be expressed as scalings of coordinates. The computation of these scalings is almost trivial.

As a side result, we get an easy proof for the fact that the group of birational automorphisms of a curve of positive genus fixing a point is always abelian. This is not a new result, as it can also be obtained as a consequence of the uniformization theorem \cite{weyl13}.

In Section \ref{SECTION: motivation and example} we detail the elliptic curve construction and motivate our interest in fixing generators. In Sections \ref{SECTION: W curves} and \ref{SECTION: TW curves} we define the Tschirnhaus-Weierstrass curves of the title and show that maps between them have a very simple form. Finally, Section \ref{SECTION: computational aspects} offers some algorithmic considerations.

\section{Motivation and example}\label{SECTION: motivation and example}

We describe in detail the construction of the equation of an elliptic curve that arises from considering the dimensions of its Riemann-Roch spaces. In what follows, $k$ is an algebraically closed field of characteristic zero, $C$ is an algebraic curve over $k$ and $P$ is a place of $C$ (for example a smooth point on the curve).

\begin{definition}\label{DEF: R-R space}
Let $C$ be an algebraic curve and $P\in C$ a place in it. For any non-negative integer $i$ we define $\mL(iP)$ as the vector space of all rational functions defined on $C$ that have a pole of order $\leq i$ at $P$ and no other poles. This is a particular instance of so called \emph{Riemann-Roch space} of $C$.
\end{definition}

\begin{example}\label{EX: elliptic curve equation from R-R spaces}
Let $C$ be elliptic. Then $\mL(0)=k$ and, by the Riemann-Roch theorem, for all $i\geq1$ and all $P\in C$, $\dim \mL(iP)=i$. We will now write bases of these vector spaces.

First $\mL(0)=\mL(P)=\langle 1\rangle$. Let $f$ be a rational function with a pole of order exactly two, then $\mL(2P)=\langle 1,f\rangle$. Similarly, let $g$ have a pole of order exactly three, then $\mL(3P)=\langle 1,f,g\rangle$. Notice that, as $k(C)$ has transcendence degree one over $k$, $f$ and $g$ are algebraically dependent.

We can write further bases without introducing more functions: we have $\mL(4P)=\langle 1,f,g,f^2\rangle$, $\mL(5P)=\langle 1,f,g,f^2,fg\rangle$, $\mL(6P)=\langle 1,f,g,f^2,fg,g^2,f^3\rangle$. But in the last space we have seven elements and dimension six, so there is a non-trivial linear combination that is equal to zero:
\begin{equation}\label{EQ: elliptic W curve}
a_1+a_2f+a_3g+a_4f^2+a_5fg+a_6g^2+a_7f^3=0, \qquad a_6,a_7\neq0.
\end{equation}

The map $P\mapsto(f(P),g(P))$ is birational from $C$ to its image in $k^2$ (see Theorem \ref{TH: W map is birational}). Thus, we have an equation of $C$ in the $(f,g)$-plane. By performing some substitutions we can eliminate certain coefficients:
\begin{itemize}
 \item By $g\leftarrow g-\frac{a_5}{2a_6} f$ we can assume that $a_5=0$.
 \item By $f\leftarrow f-\frac{a_4}{3a_7}$ we can assume that $a_4=0$.
 \item By $g\leftarrow g- a_3/(2a_6)$ we can assume that $a_3=0$.
 \end{itemize}
Note that the order is relevant since these operations change several coefficients at a time. The equation has been reduced to the classical
\begin{equation}\label{EQ: elliptic TW curve}
 a_1+a_2f+a_6g^2+a_7f^3=0, \qquad a_6,a_7\neq0.
\end{equation}
\end{example}

\begin{remark}
If we substitute $f$ or $g$ in Equation \eqref{EQ: elliptic TW curve} by scalar multiples of them, no new terms appear. But if we substitute $f$ (resp. $g$) by any combination $\alpha f+\beta$, $\beta\neq0$ (resp. $\alpha g+\beta f+\gamma$, $\beta\neq0$ or $\gamma\neq0$) then some of the terms that we eliminated appear again. In this sense, Equation \eqref{EQ: elliptic TW curve} fixes $f,g$ up to product by scalars.
\end{remark}

Our interest in these normal forms arises from the fact that curve isomorphisms have a particularly simple form. We illustrate this idea with the elliptic case.

\begin{example}\label{EX: elliptic curve isoms are scalings}
Given two elliptic curves $C_1,C_2$ and two points $P_1\in C_1,P_2\in C_2$, we want to compute all isomorphisms $\varphi:C_1\rightarrow C_2$ with
$\varphi(P_1)=P_2$. As in Example \ref{EX: elliptic curve equation from R-R spaces} we get:
\[\begin{array}{c@{}c@{}c@{}c@{}c@{\qquad}c}
a_1+&a_2f_1+&a_6g_1^2+&a_7f_1^3&=0 &\mbox{at } P_1 \\[1ex]
b_1+&b_2f_2+&b_6g_2^2+&b_7f_2^3&=0 &\mbox{at } P_2
\end{array}\]
As $\varphi$ induces linear isomorphisms $\varphi_i^*:\mL(iP_2)\rightarrow\mL(iP_1)$
for all $i$, we have that
\[\begin{array}{l@{\qquad}l}
f_2\mapsto u_1+u_2f_1, &u_2\neq0,\\[1ex]
g_2\mapsto u_3+u_4f_1+u_5g_1, &u_5\neq0.
\end{array}\]
Then we have the two equations
\[\begin{array}{l@{}l@{}l@{}l@{}l}
a_1&+a_2f_1&+a_6g_1^2&+a_7f_1^3&=0 \\[1ex]
b_1&+b_2(u_1+u_2f_1)&+b_6(u_3+u_4f_1+u_5g_1)^2&+b_7(u_1+u_2f_1)^3&=0
\end{array}\]
from which we get proportionality relations that imply $u_1,u_3,u_4=0$, that is, $f_2\mapsto u_1f_1$ and $g_2\mapsto u_5g_1$, i.e. the isomorphism is an scaling in both variables.
\end{example}

For general curves, we will have a similar result: any isomorphism from a T(schirnhaus)-W(eierstrass) curve (defined in Section \ref{SECTION: TW curves}) to another TW-curve has to be an scaling in all the variables (that is, given by a diagonal linear map). This is Theorem \ref{TH: isoms of TW curves are scalings}. The interest of such a result is in the following problem: given two curves, compute all the isomorphisms from one to the other. One solution is:

\begin{enumerate}
 \item Compute all the Weierstrass places of both curves. There are finitely many, see Definition \ref{DEF: W place} and the remark preceding it.
 \item Fix one Weierstrass place $P$ in the first curve.
 \item For each Weierstrass place $Q$ in the second curve, compute all the isomorphisms that send $P$ to $Q$.
\end{enumerate}

The last step can be done efficiently thanks to the aforementioned theorem.

\section{Weierstrass curves}\label{SECTION: W curves}

\begin{definition}\label{DEF: pole and gap numbers}
Let $C$ have genus $g\geq1$ and $P$ be a place of $C$. An integer $i\geq1$ is called a \emph{pole number} iff there exists a rational function over $C$ such that its only pole is $P$ and the order of $P$ as a pole of $f$ is $i$. Otherwise, $i$ is called a \emph{gap number}.
\end{definition}

\begin{remark}\label{REM: RR spaces grow dimension +1 or +0}
For any $i\geq1$, $\dim\mL(iP)=\dim\mL((i-1)P)+\varepsilon$, where $\varepsilon$ is 1 iff $i$ is a pole number and 0 iff $i$ is a gap number.
\end{remark}

\begin{remark}\label{REM: facts about gaps and poles}
The following are well-known facts about pole and gap numbers:
\begin{enumerate}
 \item[(i)] 1 is always a gap number: $\mL(0)=k$ always, and $L(P)=k$ for any $P$ except when $g=0$.
 \item[(ii)] By Riemann-Roch, $\dim\mL(nP)=n+1-g$ for $n\geq 2g-1$, so there are exactly $g$ gap numbers at every place by Remark \ref{REM: RR spaces grow dimension +1 or +0}.
 \item[(iii)] By the previous two items, all the places of an elliptic curve have the same gap number sequence, namely 1.
 \item[(iv)] The set of pole numbers is a semigroup with respect to addition, for the product of functions with poles at $P$ of order $i_1$ and $i_2$ is a function with a pole of order $i_1+i_2$. It is finitely generated (this is true in general for semigroups of natural numbers, in our case it is particularly easy to prove since its complement is finite).
 \item[(v)] For a generic place of a genus $g$ curve the gaps numbers are $1,2,\ldots,g$, see \cite[p. 273]{griffiths-harris78}.
\end{enumerate}
\end{remark}

The next definition is here for the sake of completeness; in what follows we do not assume that $P$ satisfies it.

\begin{definition}\label{DEF: W place}
A place in a curve of genus $g$ is called a \emph{Weierstrass place of the curve} iff its gap number sequence is not $1,\ldots,g$ (i.e. there is some gap number greater than $g$). Note that, according to the previous remark, there are no Weierstrass places if $g=1$, and places in any algebraic curve are generically non-Weierstrass.
\end{definition}

In view of the previous properties, we introduce our main definitions.

\begin{definition}
Given a place $P$ of a curve $C$ (possibly at infinity), let $R_P:=\cup_{i=0}^\infty\mL(iP)$ be the set of all rational functions with no poles other than $P$. The \emph{pole semigroup of $P$} is the set $\Gamma_P=\{n\in\Z_{\geq0} : \exists f\in R_P \mbox{ with } o_P(f)=n\}$ where $o_P(f)$ denotes the pole order of $f$ at $P$. Note that it is indeed a semigroup and that it contains 0. We will denote the elements of the minimal generating set as $d_1,\ldots,d_r$ in increasing order. When $P$ is the only place of $C$ at infinity, we will denote the semigroup simply as $\Gamma$.
\end{definition}

Note that the pole semigroup is precisely the complementary of the gap sequence.

\begin{definition}\label{DEF: W curve}
An affine curve $C\in k^n$ is called a \emph{Weierstrass curve or W-curve} when
 \begin{enumerate}
 \item it has exactly one place at infinity,
 \item the pole order of each coordinate function $x_i$ at that place is precisely $d_i$.
 \end{enumerate}
\end{definition}

\begin{definition}\label{DEF: normalised W-curve}
A W-curve is called \emph{normalised} when there exists an uniformizing parameter $t$ at infinity such that the initial coefficient in the Laurent expansion of each $x_i$ with respect to $t$ is 1.
\end{definition}

\begin{definition}
Let $C$ be a W-curve. A \emph{monomial} is a product of variables $x_1,\ldots,x_r$. The \emph{degree of a monomial} is defined as the pole order of the rational function that it represents, an element of $\Gamma$.
\end{definition}

\begin{definition}
A \emph{Tschirnhaus map or T-map} is a function $k^r\to k^r$ such that for all $i$,
\[
 x_i\mapsto x_i +\mbox{monomials of degree}<d_i.
\]
The image of a curve by a T-map is called a \emph{T-transform} of the curve.
\end{definition}

\begin{lemma}\label{LEMMA: T-transforms of W-curves are W-curves}
Every T-transform of a W-curve is a W-curve.
\end{lemma}

\begin{proof}
Clear from the properties of the degree function.
\end{proof}

\begin{lemma}\label{LEMMA: T-transforms of norm W-curves are norm W-curves}
Every T-transform of a normalised W-curve is a normalised W-curve.
\end{lemma}

\begin{proof}
By pulling back the normalising parameter for the first curve along the inverse of the T-map, one obtains a normalising parameter for the second curve.
\end{proof}

Here is a useful characterisation of normalised W-curves.

\begin{lemma}\label{LEMMA: norm curve -> same coeffs for every param}
The W-curve $C$ is normalised if and only if for any monomials $m_1,m_2$ of the same degree $d$, the function $m_1-m_2$ has pole order less than $d$.
\end{lemma}

\begin{proof}
The ``only if'' statement is clear. For the ``if'' part, let $u$ be the quotient of two monomials $m_1/m_2$ such that $\deg(m_1)-\deg(m_2)=1$; such two monomials exist because the group generated by $\Gamma$ is equal to $\Z$. Then $u$ is a uniformising parameter. We will show that it satisfies the condition of the definition of normalised W-curve.

For each $i$ consider the function $x_iu^{-d_i}=(x_im_2^{d_i})/(m_1^{d_i})$. The numerator and denominator are monomials of the same degree $d'$, so the hypothesis implies that their coefficients at degree $d'$ with respect to the parameter $u$ are equal. Thus the coefficient of $x_iu^{-d_i}$ at degree 0 is equal to 1, but then multiplying by $u^{d_i}$ we have that $x_i$ also has its coefficient of degree $d_i$ equal to 1 with respect to $u$ as required.
\end{proof}

Our goal is to find Weierstrass curves of the simplest possible form that are birational to the given curve $C$. To this end, we apply transformations that will cancel out some of their coefficients. This will be done in the next section. The last statement in this section is that for every curve, there exists a birationally equivalent W-curve.

Let $p_1,\ldots,p_r$ be the minimal generators of the semigroup of pole numbers of $P$. If $g>1$, these numbers depend on $P$. Let $f_i\in\mL(p_iP)\setminus \mL((p_i-1)P)$, $i=1,\ldots,r$. Note that one can choose each $f_i$ not only up to a constant factor, but up to a linear combination of lower order functions (for methods of computing these spaces, see \cite{hess02,hess04}). In the next section we describe a canonical form, that we call \emph{Tschirnhaus-Weierstrass curve}, with the property that the elements $f_1,\ldots,f_r$ are uniquely determined up to multiplication by constants.

\begin{theorem}\label{TH: W map is birational}
Let $C$ be a curve. Let $P$ be a place of $C$. Let $d_1,\dots,d_r$ be the minimal generating set of the pole group $\Gamma_P$. For $i=1,\ldots,r$, let $f_i\in\mL(p_iP)\setminus \mL((p_i-1)P)$. Then the rational map $\phi\colon Q\mapsto (f_1(Q),\dots,f_r(Q))$ is birational onto its image and this image is a W-curve.
\end{theorem}

\begin{proof}
We need to show that every rational function on $C$ can be expressed as a rational function in $f_1,\dots,f_r$. It is clear that every function in $R_P=\cup_{i=0}^\infty\mL(iP)$ can be expressed as a polynomial in $f_1,\dots,f_r$. Therefore it suffices to show that every nonzero rational function is a quotient of two nonzero functions in $R_P$.

Let $h$ be an arbitrary nonzero rational function. Let $D$ be its divisor of poles. Let $m$ be the degree of $D$. Let $l$ be a number such that $\mL(lP)$ has dimension bigger than $m$. We claim that there exists a function $q\in\mL(lP)$ which vanishes along $D$. The claim follows from the observation that vanishing along $D$ imposes $m$ linear conditions, and $\dim(\mL(lP))>m$. But now, $g:=h\cdot q$ has no poles other than $P$, hence both $g$ and $q$ are in $R_P$.

Finally there is precisely one place at infinity since the $f_i$ can only have poles at $P$, and the coordinate functions of the image clearly satisfy the condition of W-curve.
\end{proof}

\section{Tschirnhaus-Weierstrass curves}\label{SECTION: TW curves}

In this section we define Tschirnhaus-Weierstrass curves (abbreviated TW-curves) by imposing cancellation of certain monomials. To this end we need to define univocally a concept of coefficient. There are at least two approaches for this. One is using a local uniformizing parameter, with respect to which we can look at coefficients in Laurent expansions, but this is depentent on the parameter and it is not clear how to proceed when considering maps between curves. Instead, we will define a set of monomials (normal forms in some sense) so that every regular function is a linear combination of normal forms.

\begin{definition}\label{DEF: normal forms and coeffs}
Given a W-curve $C$, we define a set of \emph{normal forms} $N=\{n_i : i\in\Gamma\}$ where $n_i$ is an arbitrary monomial of degree $i$.

For any pair $f,c$ where $f$ is a nonzero element of $k(C)$ and $c$ is a nonnegative integer,  of degree $d$, $f$ can be expressed uniquely as a linear combination of the $n_i$ with $i\leq\deg(f)$, and we define $\coeff_C(f,c)$ to be the coefficient of the monomial $n_c$ in that linear combination. We will drop the subindex when there is no possible confusion about the curve on which this is defined.
\end{definition}

\begin{remark}
For every $d\in\Gamma$, $\coeff(n_d,d)=1$ and $\coeff(n_d,i)=0$, $i\in\Gamma,i<d$. On the other hand, for every monomial $m$ we have $\coeff(m,\deg(m))\neq0$ but not necessarily equal to 1. The normalisation introduced in Definition \ref{DEF: normalised W-curve} ensures that $\coeff(m,\deg(m))=1$.
\end{remark}

\begin{lemma}
If $C$ is a normalised W-curve then for every monomial $m$ we have $\coeff(m,\deg(m))=1$.
\end{lemma}

\begin{proof}
A consequence of Lemma \ref{LEMMA: norm curve -> same coeffs for every param}.
\end{proof}

The cancellation of monomials mentioned before takes the following form.

\begin{definition}\label{DEF: crit and TW curve wrt crit}
Given $N$, let $Crit$ be a set of triples $(m_1,m_2,k)$ such that:
\begin{itemize}
 \item $m_1$ and $m_2$ are monomials of the same degree $d$,
 \item $k\in\Z_{>0}$ such that $d-k\in\Gamma$ (in particular $k\leq d$).
\end{itemize}
A Weierstrass curve $C$ is \emph{Tschirnhaus with respect to $Crit$} or \emph{a TW-curve with respect to $Crit$} iff for every $w\in Crit$,
\[
 \cond_C(w)=0
\]
where $\cond(w)=\coeff(\coeff(m_2,d)\cdot m_1-\coeff(m_1,d)\cdot m_2,d-k)$.
\end{definition}

Note that it would have been possible to use different sets of normal forms independently for each element of $Crit$. We have not seen any use for this added flexibility, however.

\begin{definition}\label{DEF: M_crit,k and good crit}
For each $k\in\Z_{>0}$ define $Crit_k$ as the subset of elements of $Crit$ whose third component is precisely $k$ and the matrix $M_{Crit,k}$ as follows:
 \begin{itemize}
  \item The rows are indexed by the elements $w^{(i)}$ of $Crit_k$.
  \item The columns are indexed by all variables $x_j$ such that $d_j-k\in\Gamma$. Note that these do not depend on $Crit$.
  \item The $(i,j)$-th entry is equal to the exponent of $x_j$ in $m_1^{(i)}/m_2^{(i)}$.
 \end{itemize}
Then, we say that $Crit$ is \emph{good} when $M_{Crit,k}$ is square and nonsingular for every $k$.
In particular, $Crit_k$ is empty for all $k>d_r$.
\end{definition}

\begin{theorem}\label{TH: existence and uniqueness of images of NW-curves that are T-curves for good Crits}
Let $C$ be a normalised W-curve and a good $Crit$. Then there exists a T-transform of $C$ which is Tschirnhaus with respect to $Crit$. The restriction of this transformation to $C$ is unique.
\end{theorem}

\begin{proof}
Existence: we prove that for every $k$ there exists a T-transform $C_k=\phi_k(C)$ such that $\cond_{C_k}(w)=0$ for all $w\in Crit_l, l\leq k$. This is trivially true for $k=0$. Assume it true for $k-1$. We consider a map defined on $C_{k-1}$ of the form $(x_1,\ldots,x_r)\mapsto(x_1+c_1n_{d_1-k},\ldots,x_r+c_rn_{d_r-k})$ with $c_1,\ldots,c_r$ indeterminate constants. Since the coefficients corresponding to values $<k$ are preserved by this map, the image $C_k$ satisfies the required conditions up to $k-1$. Note that $C_k$ is also a normalised W-curve by Lemma \ref{LEMMA: T-transforms of norm W-curves are norm W-curves}.

Now, let $w=(m_1,m_2,k)\in Crit_k$ with $m_j=\prod x_i^{e_{i,j}}$, $j=1,2$. By Lemma \ref{LEMMA: norm curve -> same coeffs for every param}, the relation between the corresponding coefficients is precisely
\[
 \cond_{C_k}(w) = \cond_{C_{k-1}}(w) + \sum_{\substack{i=1,\ldots,r \\ d_i-k\in\Gamma}} c_i(e_{i,1}-e_{i,2}).
\]
But the entries of $M_{Crit,k}$ are precisely the coefficients of the $c_i$ for the different $w\in Crit_k$. Therefore there exists a solution $(c_1,\ldots,c_r)$ such that $\cond_{C_k}(w)=0$ for every $w\in Crit_k$ and the induction is complete.

Uniqueness: assume that $C,C'$ are TW-curves with respect to $Crit$ and $\psi\colon C\rightarrow C'$ is a T-map. We claim that $\psi$ is the identity. Assume indirectly that there exist $i_0,k_0>0$ such that the $i_0$-th component of $\psi$ is equal to $x_{i_0}+\alpha\cdot n_{d_{i_0}-k_0}+f$ with $\deg(f)<d_{i_0}-k_0$ and $\alpha\neq 0$, and let $k_0$ be as small as possible. As before, for any $w=(m_1,m_2,k_0)\in Crit_{k_0}$ with $m_j=\prod x_i^{e_{i,j}}$, $j=1,2$ we have
\[
 \cond_{C'}(w) = \cond_{C}(w) + \sum_{\substack{i=1,\ldots,r \\ d_i-k\in\Gamma}} c_i(e_{i,1}-e_{i,2}).
\]
Since $\cond_{C'}(w) = \cond_{C}(w) = 0$ and the matrix $M_{Crit,k}$ is nonsingular it follows that all the $c_i$ must be zero. However $c_{i_0}=\alpha\neq0$, contradiction.
\end{proof}

\begin{theorem}\label{TH: existence of a good Crit}
There exists a good $Crit$.
\end{theorem}

\begin{proof}
We will build a matrix $A$ with the property that we can construct all the required matrices by keeping a few well-chosen rows and columns of $A$. The construction of $Crit$ will then be straightforward.

The columns of $A$ are indexed by $x_1,\ldots,x_r$, its rows are indexed by $x_2,\ldots,x_r$. Similarly to Definition \ref{DEF: M_crit,k and good crit}, each row consists of the exponents of $x_j$ in $x_1^{d_b}/x_b^{d_1}$. That is, for $b>1$ the corresponding row is
\[
 d_b \quad 0 \quad \cdots \quad 0 \quad -d_1 \quad 0 \quad \cdots \quad 0 \\
\]

Now, fix $k>0$. As mentioned in Definition \ref{DEF: M_crit,k and good crit} the variables indexing the columns of $M_{Crit,k}$ only depend on the $d_i$ and $k$, not on $Crit$ itself. Furthermore, two possibilities arise:
\begin{itemize}
 \item if $k\neq d_1$ then the first column is not taken (the only way in which $d_1-k\in\Gamma$ is precisely $k=d_1$);
 \item if $k=d_1$ then $d_i-k\not\in\Gamma$ for any $i\geq 2$, since otherwise we would have that $d_i$ can be written in terms of smaller generators. So the only column in $M_{Crit,k}$ is the one corresponding to the first variable.
\end{itemize}
In the first case, let $x_{i_1},\ldots,x_{i_s}$ be those variables such that $d_{i_t}-k\in\Gamma$, that is, corresponding to the columns that should appear in $M_{Crit,k}$. It suffices to choose as rows those indexed by the same variables. Indeed, the result of this choice is a submatrix which is just $-d_1$ times the $s\times s$ identity matrix.

In the second case, the submatrix is a $1\times1$ matrix and we just need to choose a row with a nonzero element in the first entry, so any row suffices.

Finally, for each $k$ we have obtained pairs of monomials which make up elements $(x_1^{d_b},x_b^{d_1},k)\in Crit_k$, let $Crit$ be the union of these sets. By construction the matrices $M_{Crit,k}$ are precisely the nonsingular submatrices constructed above.
\end{proof}

An example is in order.

\begin{example}
Let $\Gamma=\langle 6,8,10,11 \rangle$. This would correspond to a Weierstrass place of gap sequence $1,2,3,4,5,7,9,13,15$ in a curve of genus 9 (we do not know of a concrete pointed curve with this property). Let the variables be $w,x,y,z$ of degrees $6,8,10,11$ respectively. We will proceed in a systematic way which can be easily implemented.

First we make a list of all possible monomials of degrees from 0 to 11. 
\[\begin{array}{c|*{5}{c}}
deg & 0 & 6 & 8 & 10 & 11 \\
\hline
mon & 1 & w & x & y & z
\end{array}\]
Since there are no two monomials of the same degree, the set $N$ is uniquely determined. This shows that every T-map must be of the form
\begin{equation}\label{EQ: T transform in example 1}
\begin{array}{l}
 w \rightarrow w + c_{1,6} \\
 x \rightarrow x + c_{2,2}\cdot w + c_{2,8} \\
 y \rightarrow y + c_{3,2}\cdot x + c_{3,4}\cdot w + c_{3,10} \\
 z \rightarrow z + c_{4,1}\cdot y + c_{4,3}\cdot x + c_{4,5}\cdot w + c_{4,11}.
\end{array}
\end{equation}
We have indexed the coefficients firstly by the variable which is transformed ($1,2,3,4$ correspond to $w,x,y,z$ resp.) and secondly by the difference between the degree of the variable and the degree of the monomial (the value of $k$ relevant to that coefficient).

The matrix $A$ from which we will extract nonsingular submatrices is
\[\begin{array}{cc}
 & \hspace{-2ex}\begin{array}{c@{\hspace{1.5em}}c@{\hspace{1.5em}}c@{\hspace{1.5em}}c} w & x & y & z \end{array} \\[1ex]
  \begin{array}{c} x \\ y \\ z \end{array} &
  \left(\begin{array}{cccc}
  8 & -6 & 0 & 0 \\
  10 & 0 & -6 & 0 \\
  11 & 0 & 0 & -6 \\
 \end{array}\right).
\end{array}\]
We need to know which columns correspond to each value of $k=1,\ldots,11$. We can read this off directly from the subindices in \eqref{EQ: T transform in example 1}, by writing down which coordinates involve which values of $k$ as their second subindex:
\[\begin{array}{c|*{11}{c}}
 k & 1 & 2 & 3 & 4 & 5 & 6 & 7 & 8 & 9 & 10 & 11 \\
 \hline
 cols & z & x,y & z & y & z & w & - & x & - & y & z
\end{array}\]
Now, as indicated in the proof of Theorem \ref{TH: existence of a good Crit}, for $k\neq d_1=6$ we choose the rows indexed by the same variables as the columns, and for $k=6$ we choose for example the first row. The following $Crit$ satisfies thus the hypothesis of Theorem \ref{TH: existence and uniqueness of images of NW-curves that are T-curves for good Crits}:
\[
 \left\{ (w^{11},z^6,1), (w^8,x^6,2), (w^{10},y^6,2), (w^{11},z^6,3), (w^{10},y^6,4), (w^{11},z^6,5), \right.
\]
\[
 \left. (w^8,x^6,6), (w^8,x^6,8), (w^{10},y^6,10), (w^{11},z^6,11) \right\}.
\]
\end{example}

\begin{remark}
A more natural definition of coefficients occurred to us but we did not manage to prove the previous results with it: instead of two monomials to make a suitable linear combination, it would be simpler to take only one monomial and express it as a linear combination of normal forms. One problem we found is that we may want to impose conditions coming from two monomials of the same degree expressed separately in terms of normal forms; this created some constraints on the choice of normal forms that we did not manage to get around of.
\end{remark}

Now we use the previous results to prove the main theorem.

\begin{definition}
A \emph{scaling} is a transformation $(x_1,\ldots,x_r)\mapsto(\lambda_1x_1,\ldots,\lambda_rx_r)$.
\end{definition}

Scalings preserve the T-property for W-curves.

\begin{lemma}
The image of a TW-curve with respect to some $Crit$ by a scaling is a TW-curve with respect to the same $Crit$.
\end{lemma}

\begin{proof}
For a given $(m_1,m_2,k)$ in $Crit$, the corresponding coefficient comes from expressing $\coeff(m_2,d)\cdot m_1-\coeff(m_1,d)\cdot m_2$ as a linear combination of normal forms. A scaling will multiply $m_1,m_2$ by certain coefficients (possibly not the same) but the new expression in $m_1,m_2$ will be a constant multiple of the old one. Since the monomials $n_i$ are also only multiplied by nonzero constants by the scaling, the property of a normal form appearing with zero coefficient is preserved.
\end{proof}

Now we can prove the main result of this article.

\begin{theorem}\label{TH: isoms of TW curves are scalings}
Let $(C_1,P_1)$ and $(C_2,P_2)$ be  pointed curves such that the pole semigroups of $P_1$ and $P_2$ are equal and they are both TW-curves with a common good $Crit$. Let $f\colon C_1\rightarrow C_2$ be an isomorphism between them sending $P_1$ to $P_2$. Then $f$ is a scaling.
\end{theorem}

\begin{proof}
It is clear that any W-curve can be normalised by a suitable scaling. Let $C'_1$ be the image by such a scaling $\sigma$, by the previous Lemma $C'_1$ is a normalised TW-curve. Consider the map $f\circ\sigma^{-1}\colon C'_1\rightarrow C_2$. Since $f$ is an isomorphism it respects the degrees so it must have the form
\[
 x_i\mapsto \alpha_ix_i +\mbox{monomials of smaller degree}
\]
which in general is not a T-map (for that, all the $\alpha_i$ must be equal to 1). However we make this into a T-map by another scaling: consider the function $g=\tau\circ f\circ\sigma^{-1}$ where $\tau(x_1,\ldots,x_r)=(\alpha_1^{-1}x_1,\ldots,\alpha_r^{-1}x_r)$. Again by the previous Lemma, the image of $g$ is a TW-curve again. Thus by Theorem \ref{TH: existence and uniqueness of images of NW-curves that are T-curves for good Crits} it must be precisely $C_1$ (it is already a TW-curve), that is, $g$ must be the identity. This implies that $f=\tau^{-1}\circ g\circ\sigma$ is a scaling.
\end{proof}

\begin{remark}
If we demand that $C_1,C_2$ are already normalized, we obtain strong conditions on the defining coefficients of $f$, see the algorithm at the end of the next section.
\end{remark}

\begin{corollary}
There is a 1-1 correspondence between the isomorphisms of two pointed curves, and the scaling isomorphisms between their TW forms.
\end{corollary}

\section{Computational aspects}\label{SECTION: computational aspects}

Assume that $C$ is the normalisation of a plane curve $C'$ given by an equation in 3 homogeneous variables. Assume that $P\in C$ is a point of $C$, given as a place of $C'$. Places of a possibly singular curve may be given by their center, a point in the plane, together with additional information that distinguishes a point in the normalisation in case the center is a singular point of $C'$. In Magma, for instance, this additional information is a pair of rational functions with only a single common zero. Then there are known algorithms to compute the semigroup $\Gamma$ at $P$ and, for each generator $d_i$, $i=1,\dots,r$, a function $f_i\in L(d_iP)$ with pole order $d_i$. This defines an isomorphism from $C$ to a Weierstrass curve as pointed out in Theorem \ref{TH: W map is birational}.

One can then define a term order in the variables of the image space and compute a Gr\"obner basis of the image, using the elimination method implemented in Magma (by default, reverse lexicographic order is chosen). The set $N$ is then the set of reduced monomials with respect to that term order, and the coefficient functions can be computed by the normal form algorithm.

The next step is to find a scaling such that the image curve is normalised. Let $m$ be a product of coordinates and their multiplicative inverses which has degree 1 -- it exists because the group generated by $\Gamma$ is equal to ${\mathbb{Z}}$. For each coordinate $x_i$, $i=1,\dots,r$, we write the product $x_i^{-1}m^{d_i}$ as a quotient of two monomials $m_i/n_i$. Let $e_i=\deg(m_i)=\deg(n_i)$ and $\alpha_i$ be the quotient of the leading coefficients of $m_i$ and $n_i$. Then we claim the scaling $(x_1,\dots,x_r)\mapsto (\alpha_1x_1,\dots,\alpha_rx_r)$ achieves normalisation. Indeed, a straightforward calculation shows that after the scaling, the equations $\coeff(m_i,e_i)=\coeff(n_i,e_i)$ hold for $i=1,\dots,r$. And the exponent vectors of $m_i/n_i$ form a ${\mathbb{Z}}$-basis for all exponent vectors of quotients of monomials that have degree 0, hence it follows that $\coeff(m,e)=\coeff(n,e)$ for all pairs of monomials $(m,n)$ such that $\deg(m)=\deg(n)=e$.

Here is a summary of the algorithm for computing the set of isomorphisms
between two given pointed curves.

\begin{algorithmic}
 \REQUIRE{Two plane curves $C_1,C_2$, given by their equation, and a place $P_i$, $i=1,2$, on each of the curves.}
 \ENSURE{The set of birational isomorphisms $C_1\to C_2$ taking $P_1$ to $P_2$.}
 \FOR{$i=1,2$}
   \STATE Compute the pole semigroup $\Gamma_i$ of $C_i$ at $P_i$;
 \ENDFOR
 \IF{$\Gamma_1\ne\Gamma_2$}
   \STATE Return the empty set;
 \ENDIF
 \STATE Let $d_1,\dots,d_r$ be the generators of the pole semigroup;
 \STATE Fix a term order on the set of variables $x_1,\dots,x_r$;
 \FOR{$i=1,2$}
   \STATE Compute pole functions $(f_{i,1},\dots,f_{i,r})$ for each generator;
   \STATE Let $\phi_i:C_i\to W_i$ defined by $(f_{i,1},\dots,f_{i,r})$;
   \STATE Compute a Gr\"obner basis $G_i$ for the ideal of $W_i$;
   \STATE Compute a scaling $\sigma_i:W_i\to N_i$ such that $N_i$ is normalized;
   \STATE Compute a T-map $\tau_i:N_i\to T_i$ such that $T_i$ is Tschirnhaus;
   \STATE Let $G'_i$ be the reduced normed Gr\"obner basis for the ideal of $T_i$;
 \ENDFOR
 \STATE Define parametric scaling $\rho:(x_1,\dots,x_r)\mapsto(\alpha^{d_1}x_1,\dots,\alpha^{d_r}x_r)$;
 \STATE Compare coefficients of $\rho^\ast(G'_2)$ and $G'_1$ and obtain a system $S$ of equations for $\alpha$;
 \STATE Let $\{\rho_1,\dots,\rho_m\}$ be the set of all instances of $\rho$ obtained by solving $S$;
 \STATE Return $\{\phi_2^{-1}\circ\sigma_2^{-1}\circ\tau_2^{-1}\circ\rho_j\circ\tau_1\circ\sigma_1\circ\phi_1 \mid j=1,\dots,m \}$.
\end{algorithmic}

\bibliographystyle{amsplain}

\end{document}